\documentclass{amsart}

\usepackage{graphics}
\usepackage{epsfig}
\usepackage{environ}
\usepackage{amsmath}
\usepackage{,amsthm,amssymb}

\usepackage{xcolor}
\usepackage{color}
\usepackage{graphicx}
\usepackage{marvosym}

\usepackage{paralist}
\usepackage{multicol}

\usepackage{pstricks,pst-text,pst-grad,pst-node,pst-3dplot,pstricks-add,pst-poly,pst-coil} 
\usepackage{pst-fun,pst-blur} 

\usepackage{hyperref}

\newtheorem{theorem}{Theorem}
\theoremstyle{plain}

\newtheorem{definition}{Definition}

\newtheorem{notation}{Notation}

\newtheorem{remark}{Remark}

\makeatletter
\DeclareFontFamily{U}{tipa}{}
\DeclareFontShape{U}{tipa}{bx}{n}{<->tipabx10}{}
\newcommand{\arc@char}{{\usefont{U}{tipa}{bx}{n}\symbol{62}}}%

\newcommand{\arc}[1]{\mathpalette\arc@arc{#1}}

\newcommand{\arc@arc}[2]{%
  \sbox0{$\m@th#1#2$}%
  \vbox{
    \hbox{\resizebox{\wd0}{\height}{\arc@char}}
    \nointerlineskip
    \box0
  }%
}
\makeatother

\usepackage[title]{appendix}

\begin{document}

\title[Advances in Desargues Affine Plane Transforms]{Progress in 
Invariant and Preserving Transforms for the Ratio of Co-Linear Points in the Desargues Affine Plane Skew Field}

\author[Orgest ZAKA]{Orgest ZAKA}
\address{Orgest ZAKA: Department of Mathematics-Informatics, Faculty of Economy and Agribusiness, Agricultural University of Tirana, Tirana, Albania}
\email{ozaka@ubt.edu.al, gertizaka@yahoo.com, ozaka@risat.org}

\author[James F. Peters]{James F. Peters}
\address{James F. PETERS: Department of Electrical \& Computer Engineering, University of Manitoba, WPG, MB, R3T 5V6, Canada  \& Department of Mathematics, Ad\.{i}yaman University, 02040 Ad\.{i}yaman, Turkey}
\thanks{The research has been supported by the Natural Sciences \& Engineering Research Council of Canada (NSERC) discovery grant 185986, Instituto Nazionale di Alta Matematica (INdAM) Francesco Severi, Gruppo Nazionale  per le Strutture Algebriche, Geometriche e Loro Applicazioni grant 9 920160 000362, n.prot U 2016/000036 and Scientific and Technological Research Council of Turkey (T\"{U}B\.{I}TAK) Scientific Human Resources Development (BIDEB) under grant no: 2221-1059B211301223.}
\email{James.Peters3@umanitoba.ca}

\dedicatory{Dedicated to Girard Desargues and M. Berger}


\begin{abstract}
This paper introduces invariant transforms that preserve the ratio of either two or three co-linear points in the Desargues affine plane skew field. 
The results given here have a clean, geometric presentation based based Desargues affine plan axiomatics and definitions with skew field properties.  The main results in this paper, are 
(1) ratio of two and three points is \emph{Invariant} under transforms: Inversion, Natural Translation, Natural dilatation, 
Mobi\"us Transform, in a line of Desargues affine plane. 
(2) parallel projection of a pair of lines in the Desargues affine plane preserves the ratio of two and three points, 
(3) translations in the Desargues affine plane preserve the ratio of 2 and 3 points and (4) dilatation in the Desargues affine plane preserve the ratio of 2 and 3 points.
\end{abstract}


\maketitle

\section{Introduction and Preliminaries}
The foundations for the study of the connections between axiomatic geometry and algebraic structures were set forth by D. Hilbert \cite{Hilbert1959geometry}. And some classic research results in this context are given, for example, by  E. Artin \cite{Artin1957GeometricAlgebra}, D.R. Huges and F.C. Piper ~\cite{HugesPiper}, H. S. M Coxeter ~\cite{CoxterIG1969}, Marcel Berger in \cite{Berger2009geometry12}, and Robin Hartshorne in \cite{Hartshorne1967Foundations}.
In the advancement of our research concerning the connections between axiomatic geometry and the association of algebraic structures in affine planes~\cite{ZakaDilauto, ZakaFilipi2016, FilipiZakaJusufi, ZakaCollineations, ZakaVertex, ZakaThesisPhd, ZakaPetersIso, ZakaPetersOrder, ZakaMohammedSF, ZakaMohammedEndo, ZakaPeters2022DyckFreeGroup}, we introduce in this paper some results on invariant transforms that preserve the ratio of co-linear points in the Desargues affine plane skew field.  

In this paper,  we advance the study of the ratio of two and three co-linear points in the Desargues affine plane. Specifically, we study some \emph{Invariant-Transforms} and investigate transforms of the ratio of two and three co-linear points that are \emph{Invariant} in the Desargues affine plane.
In addition, we consider some transforms that \emph{preserve} the ratio of two and three points in the Desargues affine plane. Earlier, we have shown that on each line on Desargues affine plane, we can construct a skew-field simply and constructively using only simple elements of elementary geometry and basic axioms of the Desargues affine plane (see \cite{ZakaFilipi2016, FilipiZakaJusufi, ZakaThesisPhd, ZakaPetersIso} ).  Results are given for the translations, parallel projections and dilatations that preserve the ratio of two and three co-linear points.
For a characterization of dilatations in general, see~\cite[vol {I}, \S 2.5.6, p. 51]{Berger2009geometry12}.  

The novelty in this paper is that we achieve our results without the use coordinates. We make use of properties enjoyed by transformations in the Desargues affine plane such as parallel projection, translations and dilatations.  We prove that these transformations preserve the ratio of two and three co-linear points.  

\textbf{Desargues Affine Plane.} Let $\mathcal{P}$ be a nonempty space, $\mathcal{L}$ a nonempty subset of $\mathcal{P}$. The elements $p$ of $\mathcal{P}$ are points and an element $\ell$ of $\mathcal{L}$ is a line. 

\begin{definition}
The incidence structure $\mathcal{A}=(\mathcal{P}, \mathcal{L},\mathcal{I})$, called affine plane, where satisfies the above axioms:
\begin{compactenum}[1$^o$]
\item For each points $\left\{P,Q\right\}\in \mathcal{P}$, there is exactly one line $\ell\in \mathcal{L}$ such that $\left\{P,Q\right\}\in \ell$.
\item For each point $P\in \mathcal{P}, \ell\in \mathcal{L}, P \not\in \ell$, there is exactly one line $\ell'\in \mathcal{L}$ such that
$P\in \ell'$ and $\ell\cap \ell' = \emptyset$\ (Playfair Parallel Axiom~\cite{Pickert1973PlayfairAxiom}).   Put another way,
if the point $P\not\in \ell$, then there is a unique line $\ell'$ on $P$ missing $\ell$~\cite{Prazmowska2004DemoMathDesparguesAxiom}.
\item There is a 3-subset of points $\left\{P,Q,R\right\}\in \mathcal{P}$, which is not a subset of any $\ell$ in the plane.   Put another way,
there exist three non-collinear points $\mathcal{P}$~\cite{Prazmowska2004DemoMathDesparguesAxiom}.
\end{compactenum}
\end{definition}

\begin{figure}[htbp]
	\centering%
		\includegraphics[width=0.85\textwidth]{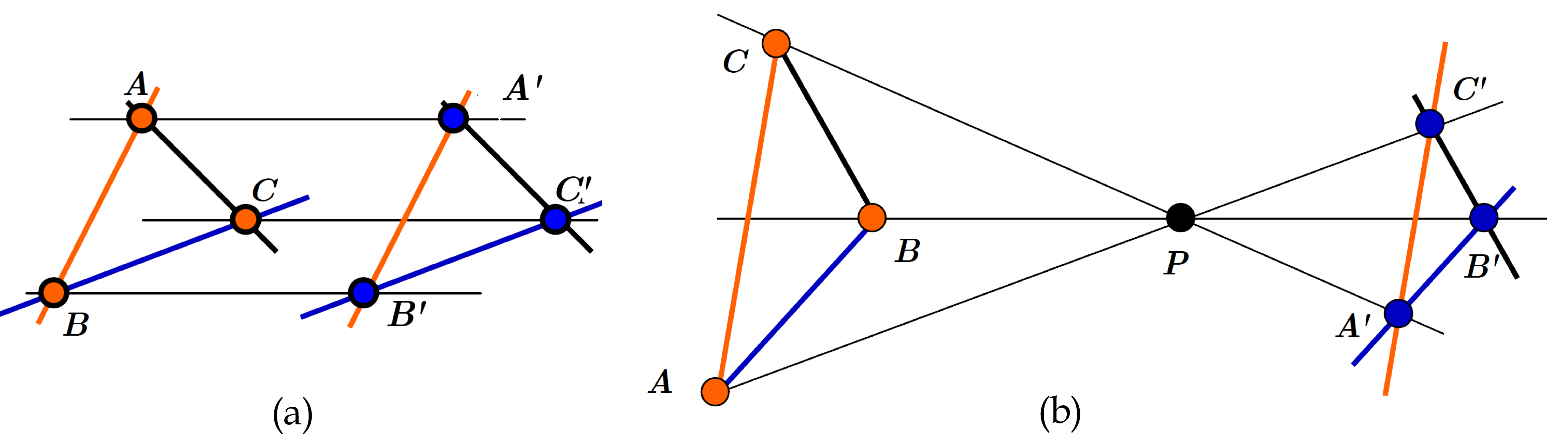}
	\caption{Desargues Axioms: (a) For parallel lines $\ell^{AA'} \parallel \ell^{BB'} \parallel \ell^{CC'}$; (b) For lines which are cutting in a single point $P$,  $\ell^{AA'} \cap \ell^{BB'} \cap \ell^{CC'}=P$.}
	\label{fig:FigureDAxiom}
\end{figure}

\emph{\bf Desargues' Axiom, circa 1630}~\cite[\S 3.9, pp. 60-61] {Kryftis2015thesis}~\cite{Szmielew1981DesarguesAxiom}.   Let $A,B,C,A',B',C'\in \mathcal{P}$ and let pairwise distinct lines  $\ell^{AA_1} , \ell^{BB'}, \ell^{CC'}, \ell^{AC}, \ell^{A'C'}\in \mathcal{L}$ such that
\begin{align*}
\ell^{AA_1} \parallel \ell^{BB'} \parallel \ell^{CC'} \ \mbox{(Fig.~\ref{fig:FigureDAxiom}(a))} &\ \mbox{\textbf{or}}\
\ell^{AA_1} \cap \ell^{BB'} \cap \ell^{CC'}=P.
 \mbox{(Fig.~\ref{fig:FigureDAxiom}(b) )}\\
 \mbox{and}\  \ell^{AB}\parallel \ell^{A'B'}\ &\ \mbox{and}\ \ell^{BC}\parallel \ell^{B'C'}.\\
A,B\in \ell^{AB}, A'B'\in \ell^{A'B'},  &\ \mbox{and}\ B,C\in \ell^{BC},  B'C'\in \ell^{B'C'}.\\
A\neq C, A'\neq C', &\ \mbox{and}\ \ell^{AB}\neq \ell^{A'B'}, \ell^{BC}\neq \ell^{B'C'}.
\end{align*}

Then $\boldsymbol{\ell^{AC}\parallel \ell^{A'C'}}$.   \qquad \textcolor{blue}{$\blacksquare$}

\noindent A {\bf Desargues affine plane} is an affine plane that satisfies Desargues' Axiom.

\begin{notation}
Three vertexes $ABC$ and $A'B'C'$, which, fulfilling the conditions of the Desargues Axiom, we call \emph{'Desarguesian'}.
\end{notation}

We, earlier, have defined, the actions: 'addition of points' and 'multiplication of points' in a line of Desargues affine planes, presented in \cite{ZakaThesisPhd, FilipiZakaJusufi, ZakaFilipi2016, ZakaVertex, ZakaPeters2022DyckFreeGroup}.

The process of construct the points $C$ for adition and multiplication of points in $\ell^{OI}-$line in affine plane, is presented in the tow algorithm form  

\begin{figure}[htbp]
\centering%
\includegraphics[width=0.80\textwidth]{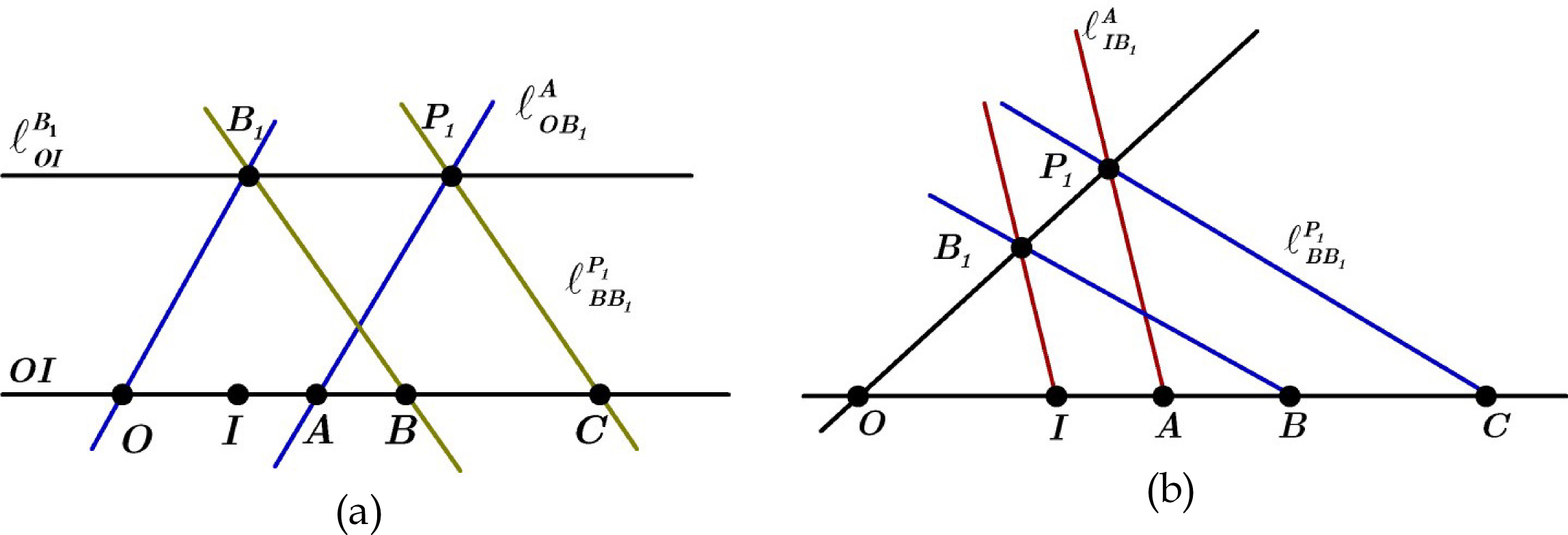}
\caption{ (a) Addition of points in a line in affine plane, 
(b) Multiplication of points in a line in affine plane}
\label{fig:FigureAdMult}
\end{figure}

\begin{multicols}{2}
\textsc{Addition Algorithm}(Fig.\ref{fig:FigureAdMult}(a))
\begin{description}
	\item[Step.1] $B_{1}\notin \ell^{OI}$
	\item[Step.2] $\ell_{OI}^{B_{1}}\cap \ell_{OB_{1}}^{A}=P_{1}$
	\item[Step.3] $\ell_{BB_{1}}^{P_{1}}\cap \ell^{OI}=C(=A+B)$
\end{description}

\textsc{Multiplication Algorithm}(Fig.\ref{fig:FigureAdMult}(b))
\begin{description}
	\item[Step.1] $B_{1}\notin \ell^{OI}$
	\item[Step.2] $\ell_{IB_{1}}^{A}\cap \ell^{OB_{1}}=P_{1}$
	\item[Step.3] $\ell_{BB_{1}}^{P_{1}}\cap \ell^{OI}=C(=A\cdot B)$
\end{description}
\end{multicols}

In \cite{ZakaThesisPhd} and \cite{FilipiZakaJusufi}, we have prove that $(\ell^{OI}, +, \cdot)$ is a skew field in Desargues affine plane, and is field (commutative skew field) in the Papus affine plane.

\begin{definition}
The parallel projection between the two lines in the Desargues affine plane, will be called, a function,
\[P_P : \ell_1 \to  \ell_2, \quad \forall A,B \in \ell_1, \quad AP_P(A) || BP_P(B)
\]
\end{definition}
It is clear that this function is a bijection between any two lines in Desargues affine planes, for this reason, it can also be thought of as isomorphism between two lines.

\begin{definition}
\cite{ZakaCollineations} Dilatation of an affine plane $\mathcal{A}=(\mathcal{P},\mathcal{L}, \mathcal{I})$, called its collineation $\delta$ such that: $\forall P\neq Q \in \mathcal{P}, \delta{(PQ)}||PQ$.
\end{definition}

\begin{definition}
\cite{ZakaCollineations} Translation of an affine plane $\mathcal{A}=(\mathcal{P},\mathcal{L}, \mathcal{I})$, called identical dilatation $id_{\mathcal{P}}$ his and every other of its dilatation, about which he affine plane has not fixed points.
\end{definition}
\textbf{Some well-known results related to translations and dilatation's in Desargues affine planes.}
\begin{itemize}
	\item 
The dilatation set $\textbf{Dil}_{\mathcal{A}}$ of affine plane $\mathcal{A}$ forms a \textbf{group} with respect
to composition $\circ$ (\cite{ZakaCollineations}).
\item The translations set $\textbf{Tr}_{\mathcal{A}}$ of affine plane $\mathcal{A}$ forms a \textbf{group} with respect
to composition $\circ$; which is a sub-group of the dilatation group $\left(\textbf{Dil}_{\mathcal{A}}, \circ\right)$ (\cite{ZakaCollineations}).
\item In a affine plane: the group $\left(\textbf{Tr}_{\mathcal{A}}, \circ\right)$ of translations is \textbf{normal
sub-group} of the group of dilatations
$\left(\textbf{Dil}_{\mathcal{A}}, \circ\right)$ (\cite{ZakaCollineations} ).
\item Every dilatation in Desargues affine plane 
	 $\mathcal{A}_{\mathcal{D}}=(\mathcal{P},\mathcal{L}, \mathcal{I})$ which leads a line in it, is an automorphism of skew-fields constructed on the same line $\ell \in \mathcal{L},$ of the plane $\mathcal{A}_{\mathcal{D}}$ (\cite{ZakaDilauto} ).
\item Every translations in Desargues affine plane 
	 $\mathcal{A}_{\mathcal{D}}=(\mathcal{P},\mathcal{L}, \mathcal{I})$ which leads a line in it, is an automorphism of skew-fields constructed on the same line $\ell \in \mathcal{L},$ of the plane $\mathcal{A}_{\mathcal{D}}$ (\cite{ZakaDilauto}).
\item Each dilatation in a Desargues affine plane, $\mathcal{A}_{\mathcal{D}}=(\mathcal{P},\mathcal{L}, \mathcal{I})$ is an isomorphism between skew-fields constructed over isomorphic lines $\ell_1, \ell_2 \in \mathcal{L}$ of that plane (\cite{ZakaPetersIso}).
\item Each translations in a Desargues affine plane, $\mathcal{A}_{\mathcal{D}}=(\mathcal{P},\mathcal{L}, \mathcal{I})$ is an isomorphism between skew-fields constructed over isomorphic lines $\ell_1, \ell_2 \in \mathcal{L}$ of that plane (\cite{ZakaPetersIso}).
\end{itemize}

%
%
%

\textbf{Ratio of two points in a line on Desargues affine plane.} 
In the paper \cite{ZakaPeters2022DyckFreeGroup}, we have done a detailed study, related to the ratio of two and three points in a line of Desargues affine plane. Below we are listing some of the results for ratio of two and three points.
  
\begin{definition} \label{ratio2points}
\cite{ZakaPeters2022DyckFreeGroup} Lets have two different points $A,B \in \ell^{OI}-$line, and $B\neq O$, in Desargues affine plane. We define as ratio of this tow points, a point $R\in \ell^{OI}$, such that,
\[R=B^{-1}A, \qquad \text{
we mark this, with,} \qquad 
R=r(A:B)=B^{-1}A
\]
\end{definition}

For a 'ratio-point' $R \in \ell^{OI}$, and for point $B\neq O$ in line $\ell^{OI}$, is a unique defined point, $A \in \ell^{OI}$, such that $R=B^{-1}A=r(A:B)$.
 
\begin{figure}[htbp]
	\centering
		\includegraphics[width=0.75\textwidth]{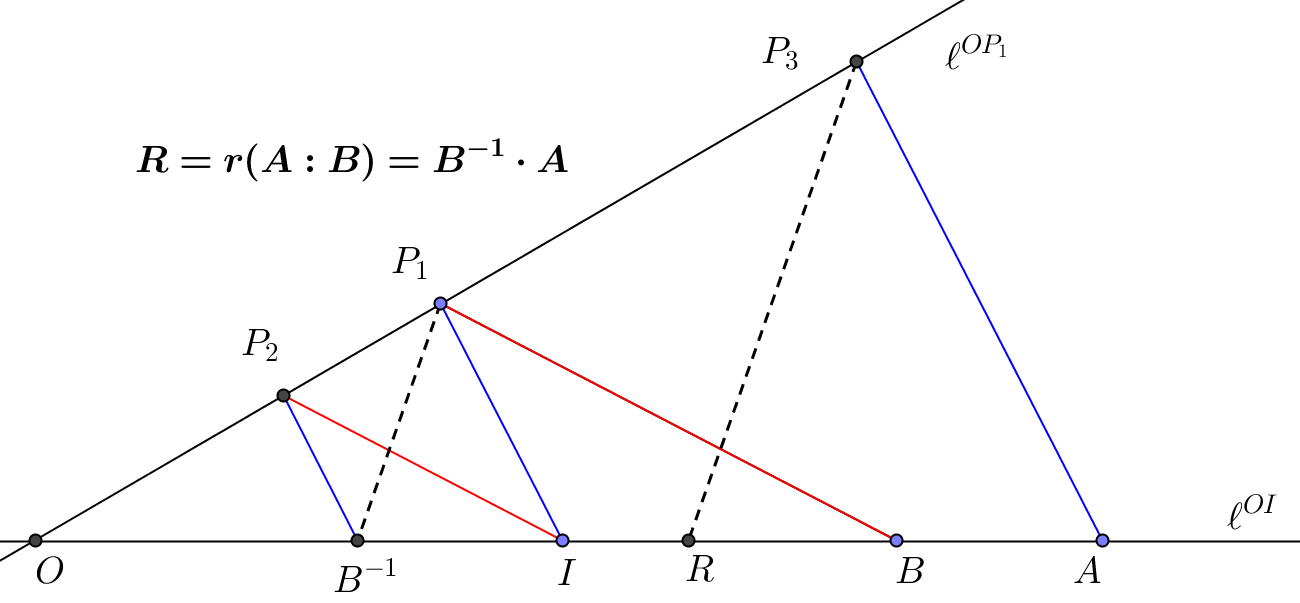}
	\caption{Ilustrate the Ratio-Point, of 2-Points in a line of Desargues affine plane $R=r(A:B)=B^{-1}A$.}
	\label{Ratio2points}
\end{figure}

\textbf{Some results for Ratio of 2-points in Desargues affine plane} (see \cite{ZakaPeters2022DyckFreeGroup}).
\begin{itemize}
\item If have two different points $A,B \in \ell^{OI}-$line, and $B\neq O$, in Desargues affine plane, then, $
r^{-1}(A:B)=r(B:A)$. 
%
\item For three collinear point $A,B,C$ and $C\neq O$, in $\ell^{OI}-$line, have, 
\[
r(A+B:C)=r(A:C)+r(B:C).
\]
\item For three collinear point $A,B,C$ and $C\neq O$, in $\ell^{OI}-$line, have,
\begin{enumerate}
	\item $r(A\cdot B:C)=r(A:C)\cdot B.$
	\item $r(A:B\cdot C)=C^{-1}r(A:C).$
\end{enumerate}
%
\item Let's have the points $A,B \in \ell^{OI}-$line where $B\neq O$.  Then have that, 
\[
r(A:B)=r(B:A) \Leftrightarrow A=B.
\]
%
\item This ratio-map, $r_{B}: \ell^{OI} \to \ell^{OI}$ is a bijection in $\ell^{OI}-$line in Desargues affine plane. 
%
\item The ratio-maps-set $\mathcal{R}_2=\{r_{B}(X)|\forall X\in \ell^{OI} \}$, for a fixed point $B$ in $\ell^{OI}-$line, forms a skew-field with 'addition and multiplication' of points. 
This, skew field $(\mathcal{R}_2, +, \cdot)$ is sub-skew field of the skew field $(\ell^{OI}, +, \cdot)$.
\end{itemize}
\textbf{Ratio of three points in a line on Desargues affine plane.} (see \cite{ZakaPeters2022DyckFreeGroup})
\begin{definition}\label{ratiodef}
If $A, B, C$ are three points on a line $\ell^{OI}$ (collinear) in Desargues affine plane, then we define their \textbf{ratio} to be a point $R \in \ell^{OI}$, such that:
\[
(B-C)\cdot R=A-C, \quad \mbox{concisely}\quad R=(B-C)^{-1}(A-C),
\]
and we mark this with  $r(A,B;C)= (B-C)^{-1}(A-C)$.
\end{definition}

\begin{figure}[htbp]
	\centering
		\includegraphics[width=0.9\textwidth]{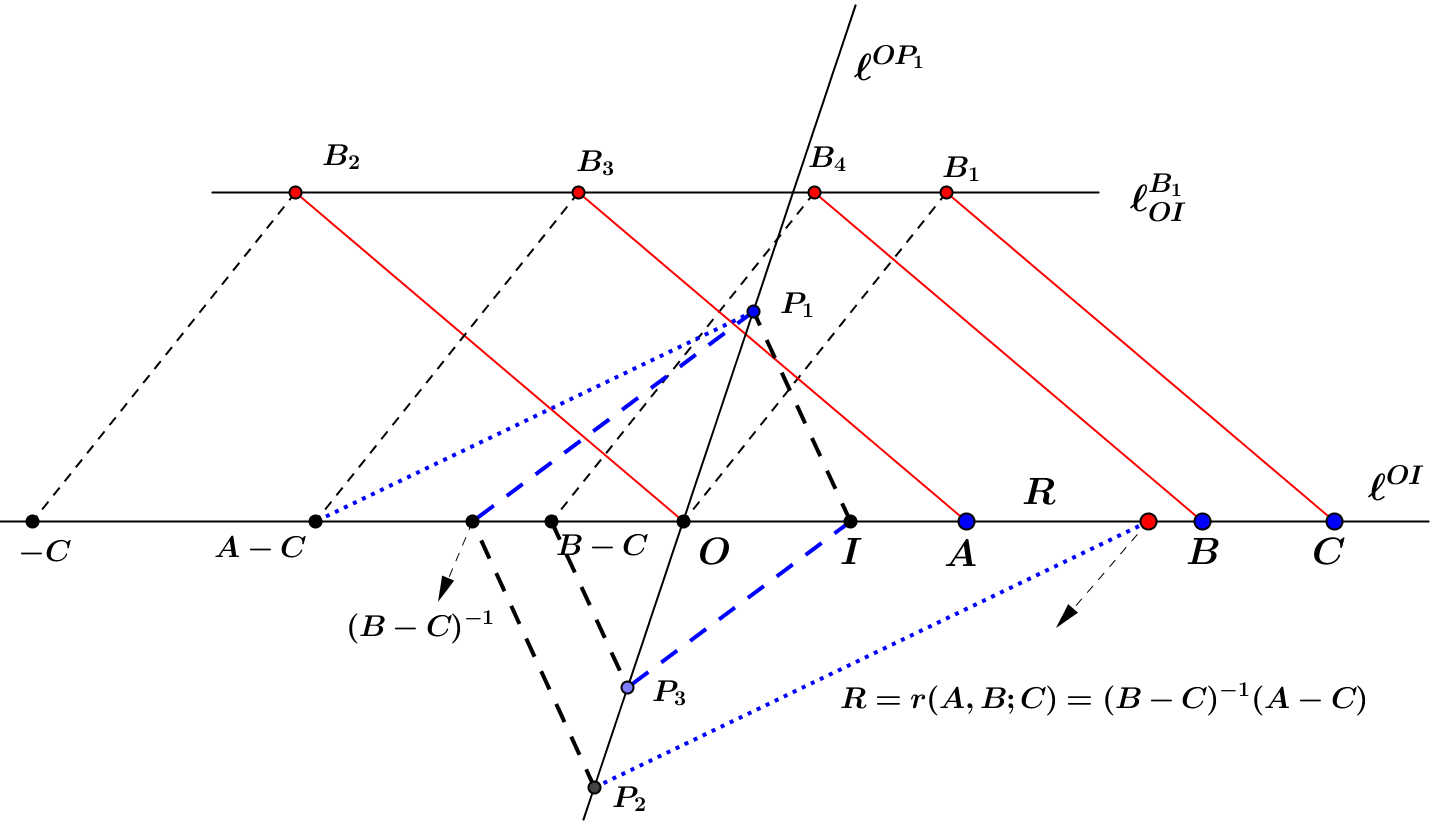}
	\caption{Ratio of 3-Points in a line of Desargues affine plane $R=r(A,B;C)$.}
		\label{ratio3points1}
\end{figure}

\textbf{Some Results for Ratio of 3-points in Desargues affine plane} (\cite{ZakaPeters2022DyckFreeGroup}).
\begin{itemize}
\item \label{reverse.ratio} For 3-points $A,B,C$ in a line $\ell^{OI}$ of Desargues affine plane, we have that,
\[
r(-A,-B;-C)=r(A,B;C).
\]
%
\item \label{inversratio} For 3-points $A,B,C$ in a line $\ell^{OI}$ in the Desargues affine plane, have
\[r^{-1}(A,B;C)=r(B,A;C).
\]
\item If  $A,B,C$, are three different points, and different from point $O$, in a line $\ell^{OI}$ on Desargues affine plane, then
\[r(A^{-1},B^{-1};C^{-1})=B[r(A,B;C)]A^{-1}.\]
%
\item In the Pappus affine plane, for three point different from point $O$, in $\ell^{OI}-$line, we have
$r(A^{-1},B^{-1};C^{-1})=r(A,B;C) \cdot r(B,A;O).$
%
%
\item This ratio-map, $r_{BC}: \ell^{OI} \to \ell^{OI}$ is a bijection in $\ell^{OI}-$line in Desargues affine plane. 
%
%
%
%
\item The ratio-maps-set $\mathcal{R}_3=\{r_{BC}(X)|\forall X\in \ell^{OI} \}$, for a different fixed points $B,C$ in $\ell^{OI}-$line, forms a skew-field with 'addition and multiplication' of points in $\ell^{OI}-$line.
This, skew field $(\mathcal{R}_3, +, \cdot)$ is sub-skew field of the skew field $(\ell^{OI}, +, \cdot)$.
\end{itemize}

\section{Invariant transforms for Ratio-2-points and Ratio-3-points}

In this section we will see some transformations, for which the ratio of 2-points and ratio of 3-points, are invariant under their action. We define these transformations first,

\begin{definition}\label{inversion}
Inversion of points in $\ell^{OI}-$line, called the map 
\[ j_P:\ell^{OI} \to \ell^{OI}, \]
which satisfies the condition, 
\[ \forall A \in \ell^{OI} \quad j_P(A)=P \cdot A.\]
\end{definition}

\begin{definition}\label{nTr}
A natural translation with point $P$, of points in $\ell^{OI}-$line, called the map 
\[ \varphi_P:\ell^{OI} \to \ell^{OI}, \]
for a fixed $P \in \ell^{OI}$ which satisfies the condition, 
\[ \forall A \in \ell^{OI} \quad \varphi_P(A)=P+A.\]
\end{definition}

\begin{definition}\label{nDil}
A natural dilatation of points in $\ell^{OI}-$line, called the map 
\[ \delta_n:\ell^{OI} \to \ell^{OI}, \]
for a fixed natural number $n \in \mathbb(N)$ which satisfies the condition, 
\[ \forall A \in \ell^{OI} \quad \delta_n(A)=nA=\underbrace{A+A+\cdots +A}_{n-times}.\]
\end{definition}

\begin{definition} \label{mobius.transform.def}
\begin{description}
	\item[(a)] For \emph{Ratio of 2-points:} Lets have an fixed points $B \in \ell^{OI}$, which are different from point $O$. Mobi\"us transform for ratio or two points in $\ell^{OI}-$line, we called the map, 
\[
m: \ell^{OI} \to \ell^{OI},
\]
which satisfies the condition, 
\[
\forall X \in \ell^{OI}, \quad m(X)=r(X:B).
\]
	
	\item[(b)]  For \emph{Ratio of 3-points:}
Lets have three fixed points $B,C \in \ell^{OI}.$ Mobi\"us transform for ratio, we called the map, 
\[
\mu: \ell^{OI} \to \ell^{OI},
\]
which satisfies the condition, 
\[
\forall X \in \ell^{OI}, \quad \mu(X)=r(X,B;C).
\]
\end{description}
\end{definition}

\begin{theorem}
Ratio of 2-points is invariant under the natural dilatation with a fixet $n\in \mathbb{N}$.
\end{theorem}
\proof
For ratio definition \ref{ratio2points} we have that, $r(A:B)=B^{-1}A $,  and for natural dilatation definition \ref{nTr}, we have,
$\delta_n(A)=nA, \quad\text{and}\quad \delta_n(B)=nB,\forall A,B \in \ell^{OI}$, so,
\[
\begin{aligned}
r[\delta_n(A):\delta_n(B)]
&=r(nA:nB)\\
&=(nB)^{-1}(nA)\\
&=B^{-1}n^{-1}nA\\
&=B^{-1}A\\
&=r(A:B).
\end{aligned}
\]
Hence  $r[\delta_n(A):\delta_n(B)]=r(A:B)$.
\qed

\begin{theorem}
Ratio of 2-points is invariant under inversion with a given point $P\in \ell^{OI}$.
\end{theorem}
\proof
For ratio definition \ref{ratio2points} we have that, $r(A:B)=B^{-1}A$, and for inversion definition \ref{inversion} with a point $P$, we have, $j_P(A)=PA, \forall A \in \ell^{OI}$, so,
\[
\begin{aligned}
r[j_P(A):j_P(B)]
&=r(PA:PB)\\
&=(PB)^{-1}(PA)\\
&=B^{-1}P^{-1}PA\\
&=B^{-1}(P^{-1}P)A\\
&=B^{-1}IA\\
&=B^{-1}A\\
&=r(A:B).
\end{aligned}
\]
Hence, $r[j_P(A):j_P(B)]=r(A:B).$
\qed

\begin{theorem}
Ratio of 2-points is invariant under Mobi\"us transform.
\end{theorem}
\proof
For Mobi\"us transform definition \ref{mobius.transform.def} we have $m(X)=r(X:B)=B^{-1}X$, so, for ratio of 2-point, the points $m(A),m(B) \in \ell^{OI}$ first, we calculate, this point, according to following the definition of $m-$map, and we have
\begin{itemize}
	\item $m(A)=r(A:B)$,
	\item $m(B)=r(B:B)$, so $	m(B)=B^{-1}B=I.$
\end{itemize}
Now, calculate,
\[
\begin{aligned}
r[m(A):m(B)] &=r(\mu(A):I)\\
&=I^{-1}[m(A)] \\
&=I\cdot m(A) \\
&=m(A)\\
&=r(A:B).
\end{aligned}
\]
Hence, $r[m(A):m(B)]=r(A:B).$
\qed

\begin{theorem}
Ratio of 3-points, is invariant under the natural translation with a point $P$.
\end{theorem}
\proof
For ratio definition \ref{ratio2points} we have that,
$r(A,B;C)=(B-C)^{-1}(A-C)$, 
and for natural translation with a point $P$ definition \ref{nTr}, we have, $\varphi_P(A)=P+A, \forall A \in \ell^{OI}$, so, for ratio we have that,
\[
\begin{aligned}
r[\varphi_P(A),\varphi_P(B);\varphi_P(C)]
&=r(A+P,B+P;C+P)\\
&=([B+P]-[C+P])^{-1}([A+P]-[C+P])\\
&=(B+P-C-P)^{-1}(A+P-C-P)\\
&=(B-C)^{-1}(A-C)\\
&=r(A,B;C).
\end{aligned}
\]
Hence, $r[\varphi_P(A),\varphi_P(B);\varphi_P(C)]=r(A,B;C).$
\qed

\begin{theorem}
Ratio of 3-points is invariant under the natural dilatation with a fixet $n\in \mathbb{N}$.
\end{theorem}
\proof
For ratio definition \ref{ratiodef} we have that,
$r(A,B;C)=(B-C)^{-1}(A-C)$, 
and for natural dilatation definition \ref{nDil}, we have,
$
\delta_n(A)=nA, \forall A \in \ell^{OI}
$,
so,
\[
\begin{aligned}
r[\delta_n(A),\delta_n(B);\delta_n(C)]
&=r(nA,nB;nC)\\
&=(nB-nC)^{-1}(nA-nC)\\
&=(n[B-C])^{-1}(n[A-C])\\
&=(B-C)^{-1}n^{-1}n(A-C)\\
&=(B-C)^{-1}(A-C)\\
&=r(A,B;C).
\end{aligned}
\]
Hence, $r[\delta_n(A),\delta_n(B);\delta_n(C)]=r(A,B;C).$
\qed

\begin{theorem}
Ratio of 3-points is invariant under Inversion with a given point $P\in \ell^{OI}$.
\end{theorem}
\proof
For ratio definition \ref{ratiodef} we have that,
$r(A,B;C)=(B-C)^{-1}(A-C)$, 
and for inversion definition \ref{inversion} with a point $P$, we have, $
j_V(A)=PA, \forall A \in \ell^{OI}
$, so,
\[
\begin{aligned}
r[j_P(A),j_P(B);j_P(C)]
&=r(PA,PB;PC)\\
&=(PB-PC)^{-1}(PA-PC)\\
&=(P[B-C])^{-1}(P[A-C])\\
&=(B-C)^{-1}P^{-1}P(A-C)\\
&=(B-C)^{-1}I(A-C)\\
&=(B-C)^{-1}I(A-C)\\
&=r(A,B;C).
\end{aligned}
\]
Hence, $r[j_P(A),j_P(B);j_P(C)]=r(A,B;C).$
\qed

\begin{theorem}
Ratios of 3-points, is invariant under Mobi\"us transform.
\end{theorem}
\proof
For Mobi\"us transform definition \ref{mobius.transform.def} we have 
\[\mu(X)=r(X,B;C)=[(B-C)^{-1}(X-C)]
\]
so, for cross-ratio of the points $\mu(A),\mu(B),\mu(C) \in \ell^{OI}$ first, we calculate, this point, according to following the definition of $\mu-$map, and we have
\begin{itemize}
	\item $\mu(A)=r(A,B;C)$,
	\item $\mu(B)=r(B,B;C) \Rightarrow 	\mu(B)=[(B-C)^{-1}(B-C)]=I \Rightarrow \mu(B)=I.$
	\item $\mu(C)=r(C,B;C) \Rightarrow \mu(C)=(B-C)^{-1}(C-C)=(B-C)^{-1}O=O \Rightarrow \mu(C)=O.$
\end{itemize}
Now, calculate,
\[
\begin{aligned}
r[\mu(A),\mu(B);\mu(C)]
&=r(\mu(A),I;O)\\
&=(I-O)^{-1}(\mu(A)-O) \\
&=(I)^{-1}(\mu(A)) \\
&=\mu(A)\\
&=r(A,B;C).
\end{aligned}
\]
Hence, $r[\mu(A),\mu(B);\mu(C)]=r(A,B;C).$
\qed
\begin{theorem}
For natural translation with a given point $P$, we have that the ratio of two points $A$ and $B$, is equal with ratio of three points $\varphi_P(A)$,$\varphi_P(B)$, $P$, so
\[
r(A:B)=r\left( \varphi_P(A),\varphi_P(B);P \right)
\]
\end{theorem}
\proof
For ratio definition \ref{ratio2points} we have that $r(A:B)=B^{-1}A$, in this equation we substitute  $B=[B+P]-P $ and  $A=[A+P]-P$, and have,
\[
\begin{aligned}
r(A:B)&=r\left( [A+P]-P : [B+P]-P \right)\\
&=\left([B+P]-P\right)^{-1} \left([A+P]-P \right)\\
&=\left(\varphi_P(B)-P\right)^{-1} \left(\varphi_P(A)-P \right)\\
&=r\left(\varphi_P(A),\varphi_P(B);P\right).
\end{aligned}
\]
because from natural translation definition \ref{nTr}, we have that $B+P=\varphi(B)$ and $A+P=\varphi(A)$ and we consider the ratio of three points definition. Hence,
\[
r(A:B)=r\left(\varphi_P(A),\varphi_P(B);P\right).
\]
\qed

\section{Transforms which Preserving Ratio}
In this section we prove that the parallel projection, Translations and dilatation of a line in itself or in isomorphic line in Desargues affine plane, \emph{preserving:} the ratio of of two and three points. The geometrical interpretations, even in the Euclidean view, are quite beautiful in the above theorems, regardless of the rather rendered figures. This is also the reason that we are giving the proofs in the algebraic sense. So we will always have in mind the close connection of skew field and a line in Desargues affine plane, and the properties of parallel projection, translations and dilatations. To prove our results, we refer to the parallel projection, translations and  dilatation's properties, which are studied in the papers \cite{ZakaThesisPhd, ZakaCollineations, ZakaDilauto, ZakaPetersIso, ZakaPetersOrder, ZakaVertex} and our achieved results, for, ratio of two and ratio of three points in \cite{ZakaPeters2022DyckFreeGroup}.

\begin{figure}[htbp]
	\centering
		\includegraphics[width=0.80\textwidth]{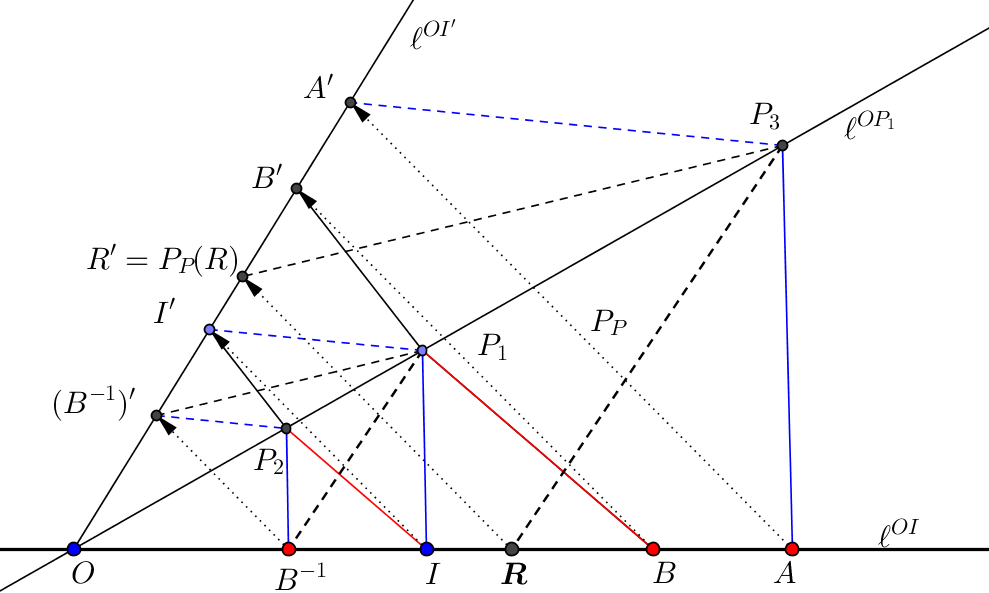}
	\caption{Ratio of 2-Points into parallel projection in Desargues affine plane.}
	\label{Ratio2points.pp}
\end{figure}

\begin{theorem}
The parallel projection between the two lines $\ell_1$ and $\ell_2$ in Desargues affine plane, \textbf{preserving the ratio} of $2-$points, 
\[
P_P(r(A:B))=r(P_P(A):P_P(B))
\]
\end{theorem}
\proof
If $\ell_1 || \ell_2$, we have that the parallel projection is a translation, and have true this theorem.\\
If lines $\ell_1$ and $\ell_2$  they are not parallel (so, they are cut at a single point), we have $A,B \in \ell_1$ and $P_P(A), P_P(B) \in \ell_2$. It is easily proved, with the help of Desargues' configuration (with Desargues axiom and affine plane axioms), that $P_P(r(A:B))=r(P_P(A):P_P(B))$. In figure \ref{Ratio2points.pp}, we have market $\ell_1=\ell^{OI}$, $\ell_2=\ell^{OI'}$, $I'=P_P(I)$, $A'=P_P(A)$, $B'=P_P(A)$, $(A^{-1})'=P_P(A^{-1})$, $(B^{-1})'=P_P(B^{-1})$, $R'=P_P(R)$.

If we consider, two 'three-vertex' $B^{-1}P_1(B^{-1})'$ and $RP_3P_P(R)$ in Fig.\ref{Ratio2points.pp}, they are Desarguesian, for this reason, and for the results for construction of points $R'=(B^{-1})'A'$ we have that,
\[
R'=P_P(R)
\]
because, $P_3R' || P_3P_P(R) || P_1P_P(B^{-1})$, and parallelism is equivalence relation. So the points, $R'$ and $P_P(R)$ are the cutting points of lines $\ell^{OI'}$ and $\ell^{P_3}_{P_1(B^{-1})'}$, but this is a single point, for this, we have that $R'=P_P(R).$ 
\qed

\begin{theorem}
Translations in Desargues affine plane, \textbf{preserving the ratio} of $2-$points $A,B$ in a line $\ell^{OI}$ of this plane, 
\[
\varphi(r(A:B))=r(\varphi(A):\varphi(B))
\]
\end{theorem}
\proof
Lets have a translation $\varphi: \ell^{OI} \to \ell^{O'I'}$, we know that translation preserves parallelism, therefore $\ell^{OI} \parallel \ell^{O'I'}$.

We mark with: $O'=\varphi(O)$, $I'=\varphi(I)$, $A'=\varphi(A)$, $B'=\varphi(B)$. Also, translation $\varphi$, we also apply it to auxiliary points, and mark $P_1'=\varphi(P_1)$, $P_2'=\varphi(P_2)$, $P_3'=\varphi(P_3)$. It is easily shown from the properties of the construction of the inverse point and from the Desargues condition, that $\varphi(B^{-1})=[\varphi(B)]^{-1}$. 

\begin{figure}[htbp]
	\centering
		\includegraphics[width=0.80\textwidth]{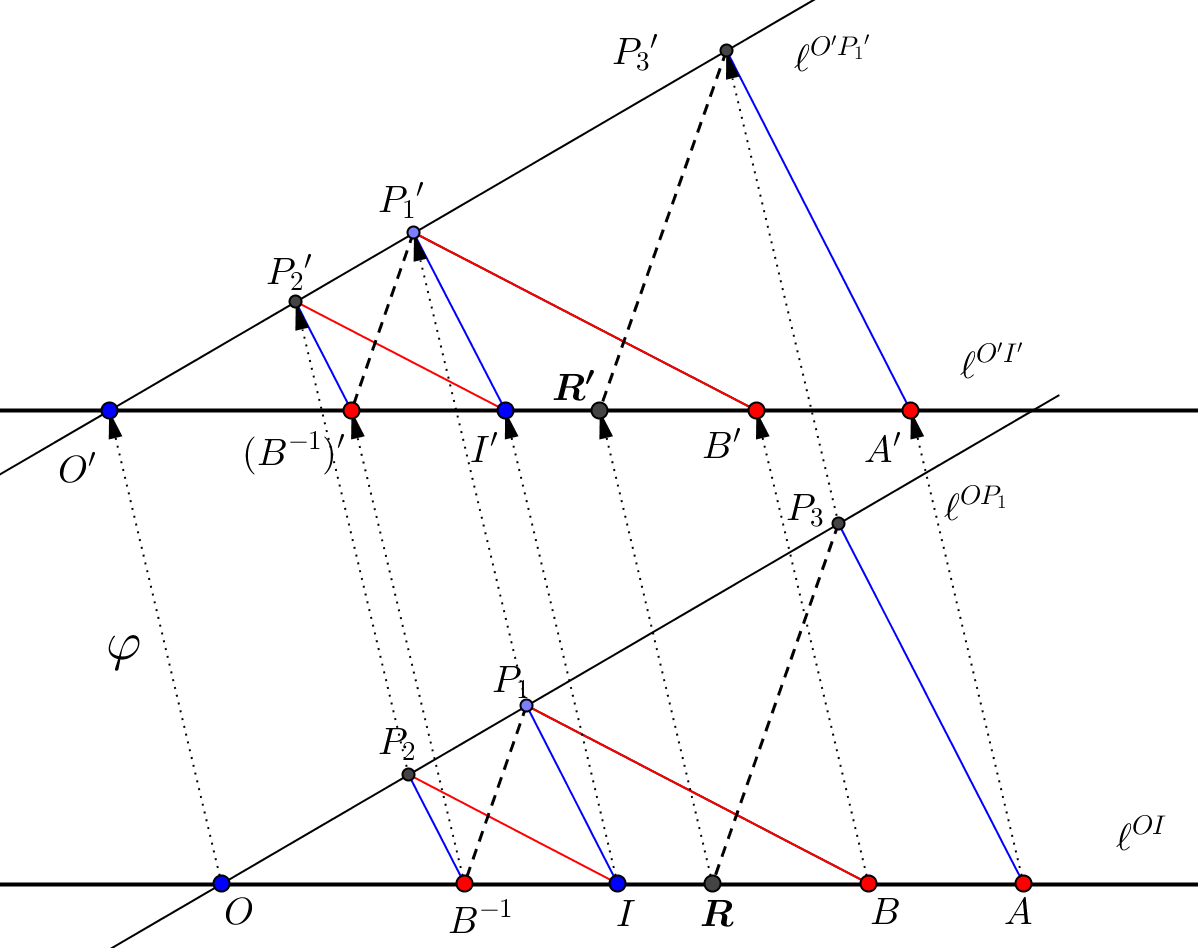}
	\caption{Ratio of 2-Points into translation in Desargues affine plane.}
	\label{fig:Tr-ratio}
\end{figure}

If we consider, two 'three-vertex' $AP_3R$ and $A'P_3'\varphi(R)$, in Fig.\ref{fig:Tr-ratio}, we see that they are \emph{Desarguesian}, for this reason, and for the results for construction of points $R'=(B^{-1})'A'$ we have that,
\[
R'=\varphi(R)
\]
because, $P_3R||P_1B^{-1}||P_1'(B^{-1})'=P_1'(B')^{-1} || P_3\varphi(R) || P_3R'$, and parallelism is equivalence relation. So the points, $R'$ and $\varphi(R)$ are the cutting points of lines $\ell^{O'I'}$ and $\ell^{P_3'}_{P_1'(B^{-1})'}$, but this is a single point, for this, we have that $R'=\varphi(R).$ 

The case when, $\varphi: \ell^{OI} \to \ell^{OI}$ (there are traces $\ell^{OI}$), $\varphi$ can be seen as a composition of two translations $\varphi_1,\varphi_2$, with different trace of the line $\ell^{OI}$.
\qed

\begin{theorem}
dilatation $\delta$ with fixed point in the same line $\ell^{OI}$ or with fixed point out of line, of points $A,B$, preserve the ratio of this points,
\[
\delta(r(A:B))=r(\delta(A):\delta(B))
\]
\end{theorem}
\proof
Lets have firstly, a dilatation with an fixed point $V \notin \ell^{OI}$ in Desargues affine plane, which $\delta: \ell^{OI} \to \ell^{O'I'}$ we know that dilatation preserves parallelism, therefore $\ell^{OI} \parallel \ell^{O'I'}$.

\begin{figure}[htbp]
	\centering
		\includegraphics[width=0.80\textwidth]{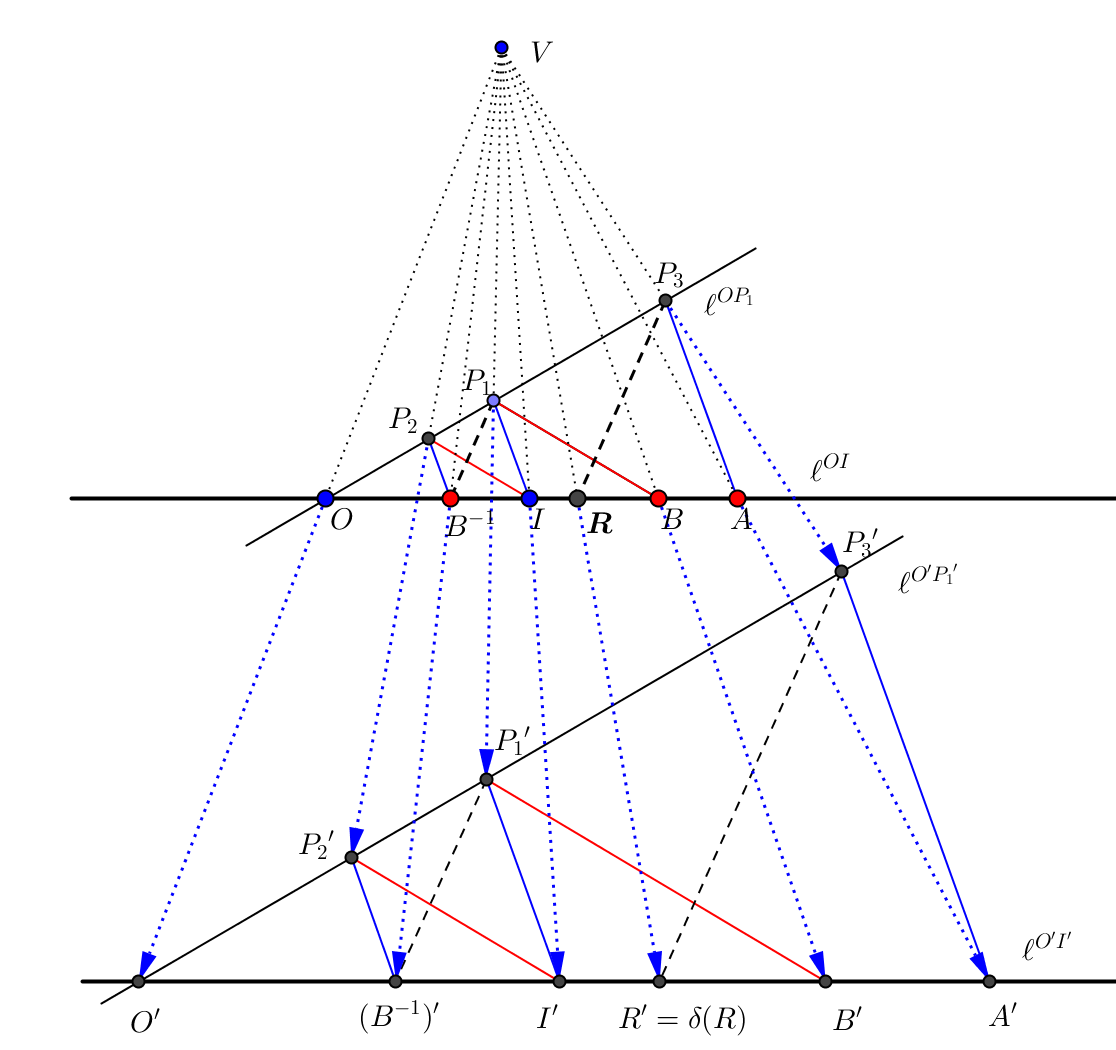}
	\caption{Ratio of 2-Points into translation in Desargues affine plane.}
	\label{Ratio2points.dil}
\end{figure}

We mark with: $O'=\delta(O)$, $I'=\delta(I)$, $A'=\delta(A)$, $B'=\delta(B)$, $\delta(R)$ Fig.\ref{Ratio2points.dil}. Also, dilatation $\delta$, we also apply it to auxiliary points, and mark $P_1'=\delta(P_1)$, $P_2'=\delta(P_2)$, $P_3'=\delta(P_3)$. It is easily shown from the properties of the construction of the inverse point and from the Desargu condition, that $\delta(B^{-1})=[\delta(B)]^{-1}=(B^{-1})'$. 

If we consider, two 'three-vertex' $AP_3R$ and $A'P_3'\delta(R)$, in Fig.\ref{Ratio2points.dil}, we see that they are Desarguesian, for this reason, and for the results for construction of points $R'=(B^{-1})'A'$ we have that,
\[
R'=\delta(R)
\]
because, $P_3R||P_1B^{-1}||P_1'(B^{-1})'=P_1'(B')^{-1} || P_3\delta(R) || P_3R'$, and parallelism is equivalence relation. So the points, $R'$ and $\delta(R)$ are the cutting points of lines $\ell^{O'I'}$ and $\ell^{P_3'}_{P_1'(B^{-1})'}$ (the line which passes from point $P_3'$ and is parallel with line $\ell^{P_1'(B^{-1})'}$), but this is a single point, for this, we have that $R'=\delta(R).$ 
\[ R'=\ell^{OI}\cap \ell^{P_3'}_{P_1'(B^{-1})'} \quad \text{and}\quad \delta(R)=\ell^{OI}\cap \ell^{P_3'}_{P_1'(B^{-1})'},
\]
so,
\[ R'=\delta(R). \]
Hence,
\[\delta \left[ r(A:B) \right]=r \left( \delta(A):\delta(B) \right).
\]

Now let's look at the case where, $\delta: \ell^{OI} \to \ell^{OI}$ (is case where the fixed point $V \in \ell^{OI}$). For simplicity of interpretation, we are taking the fixed point $V=O$, (the case where the point $V\neq O$ will be discussed at the end of the proof.

\begin{figure}[htbp]
	\centering
		\includegraphics[width=0.95\textwidth]{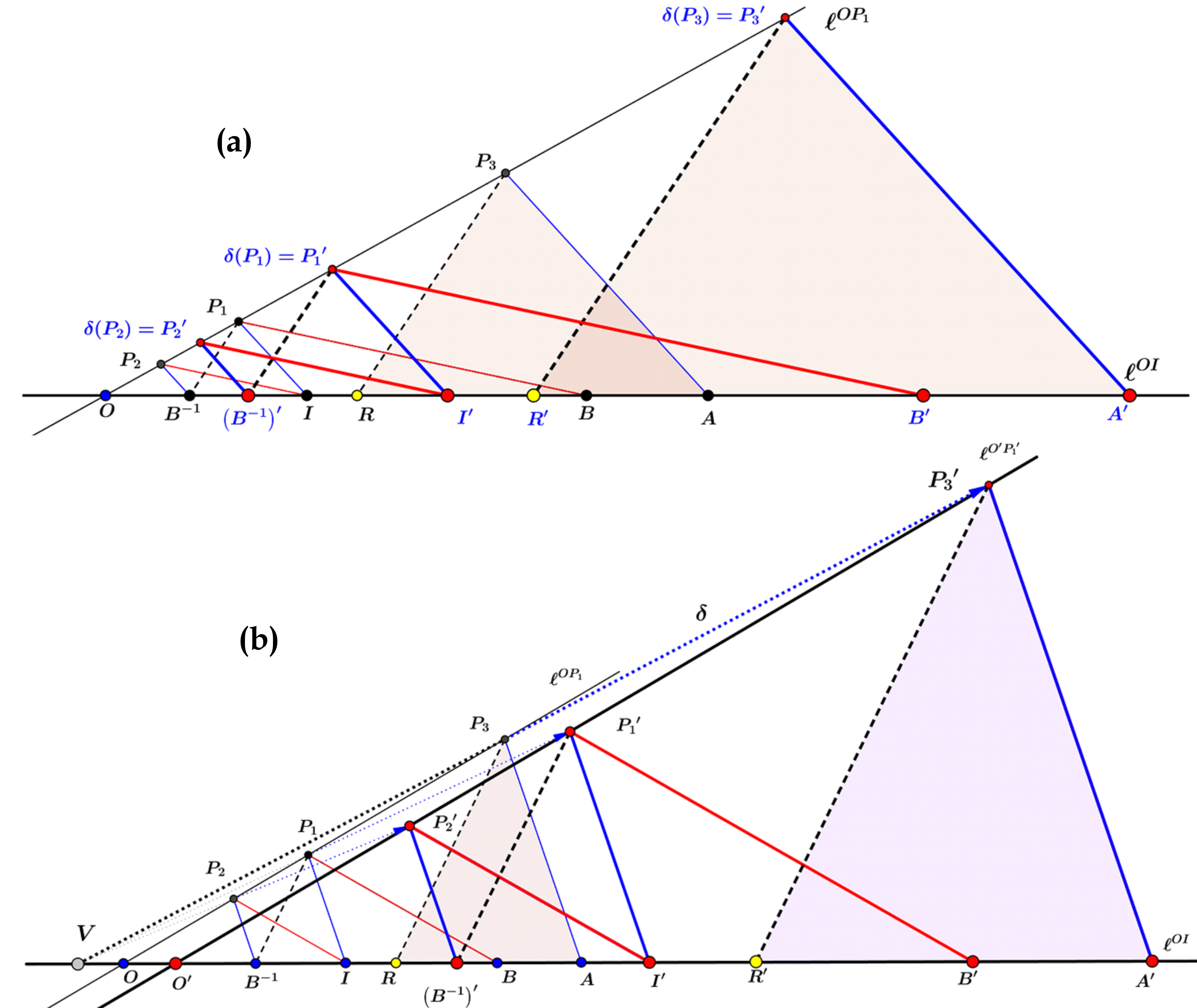}
 	\caption{Ratio of 2-Points into dilatation with fixed point in $\ell^{OI}-$line, in Desargues affine plane; (a) case where $V=O$ and (b) case where $V\neq O$.}
	\label{Ratio2points.dil0}
\end{figure}

As in the first case, we mark with: $O'=\delta(O)=O$, $I'=\delta(I)$, $A'=\delta(A)$, $B'=\delta(B)$, $\delta(R)$, see 
Fig.\ref{Ratio2points.dil0} (a),
(we are making these notes to simplify the symbolism a bit). Also, dilatation $\delta$, we also apply it to auxiliary points, and mark $P_1'=\delta(P_1)$, $P_2'=\delta(P_2)$, $P_3'=\delta(P_3)$. It is easily shown from the properties of the construction of the inverse point and from the Desargu condition, that $\delta(B^{-1})=[\delta(B)]^{-1}=(B^{-1})'$. 

If we consider, two 'three-vertex' $AP_3R$ and $A'P_3'\delta(R)$, in 
Fig.\ref{Ratio2points.dil0} (a), we see that they are Desarguesian, for this reason, and for the results for construction of points $R'=(B^{-1})'A'$ we have that,
\[
R'=\delta(R)
\]
because, $P_3\delta(R) ||  P_3R||P_1B^{-1}||P_1'(B^{-1})'=P_1'(B')^{-1}  || P_3R'$, and parallelism is equivalence relation, hence
\[
P_3R' \parallel P_3\delta(R)
\]
So the points, $R'$ and $\delta(R)$ are the cutting points of lines $\ell^{OI}$ and $\ell^{P_3'}_{P_1'(B^{-1})'}$ (the line which passes from point $P_3'$ and is parallel with line $\ell^{P_1'(B^{-1})'}$), but this is a single point, for this, we have that $R'=\delta(R).$ 
\[ R'=\ell^{OI}\cap \ell^{P_3'}_{P_1'(B^{-1})'} \quad \text{and}\quad \delta(R)=\ell^{OI}\cap \ell^{P_3'}_{P_1'(B^{-1})'},
\]
so,
\[ R'=\delta(R). \]
Hence,
\[\delta \left[ r(A:B) \right]=r \left( \delta(A):\delta(B) \right).
\]

In the same way, the case where the fixed point $V\notin O$ will be proven. The proof can be done directly see Fig.\ref{Ratio2points.dil0} (b), or this dilatation can be expressed, as i composition of an translation which con point $V$ to point $O$.
\qed

\begin{theorem}
Translations in Desargues affine plane, \textbf{preserving the ratio} of $3-$points $A,B,C$ in a line $\ell^{OI}$ of this plane, 
\[
\varphi(r(A,B;C))=r(\varphi(A), \varphi(B);\varphi(C))
\]
\end{theorem}
\proof
Lets have a translation with trace different from the $\ell^{OI}-$line, so $\varphi: \ell^{OI} \to \ell^{O'I'}$, we know that translation preserves parallelism, therefore $\ell^{OI} \parallel \ell^{O'I'}$.

\begin{figure}[htbp]
	\centering
		\includegraphics[width=0.95\textwidth]{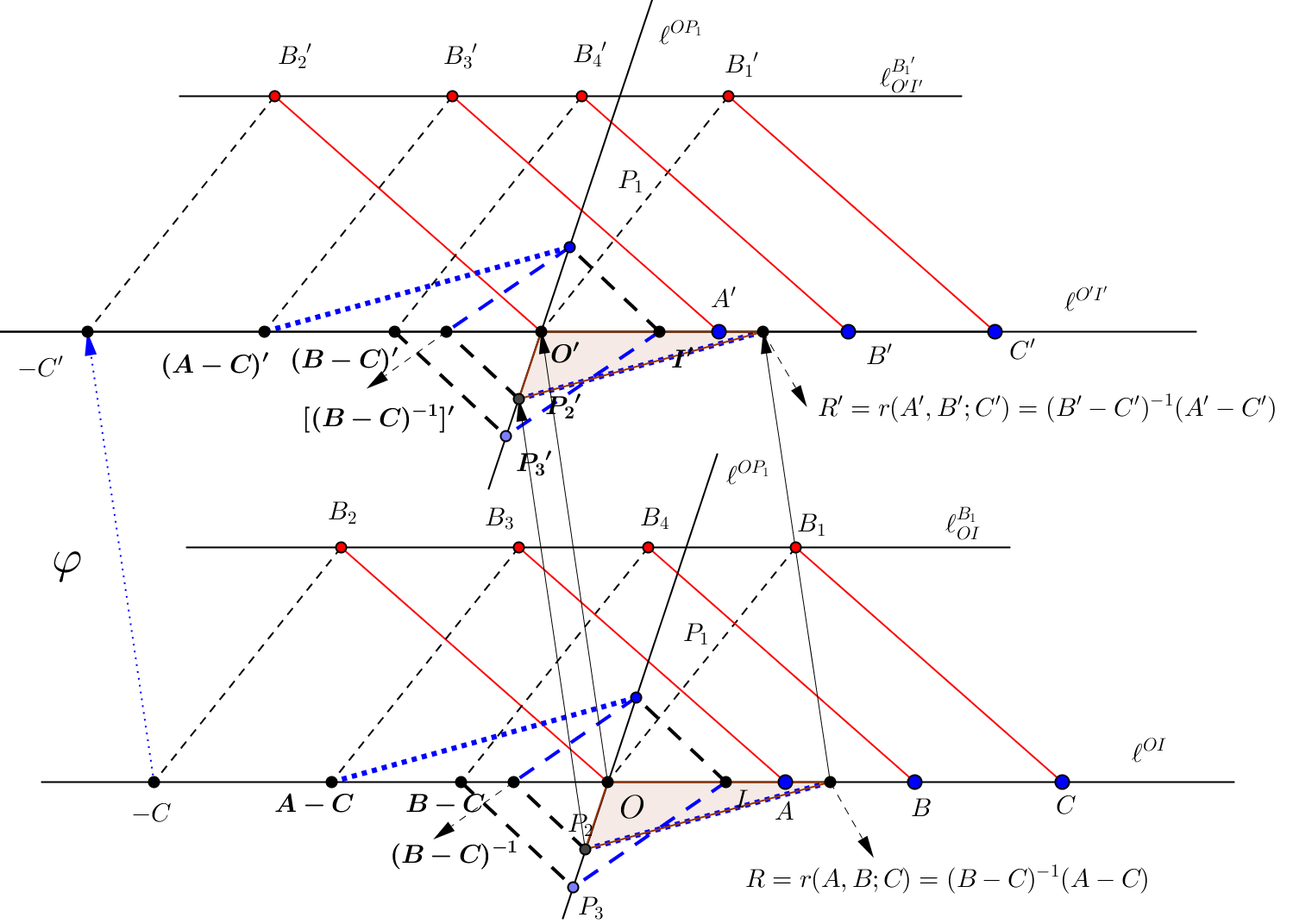}
	\caption{Ratio of 3-Points into translation in Desargues affine plane.}
	\label{Ratio3points.tr}
\end{figure}

We mark with: $O'=\varphi(O)$, $I'=\varphi(I)$, $A'=\varphi(A)$, $B'=\varphi(B)$, $C'=\varphi(C)$. Also, translation $\varphi$, we also apply it to auxiliary points, and mark $B_1'=\varphi(B_1)$, $B_2'=\varphi(B_2)$, $B_3'=\varphi(B_3)$, $B_4'=\varphi(B_4)$ and  $P_1'=\varphi(P_1)$, $P_2'=\varphi(P_2)$, $P_3'=\varphi(P_3)$. From results in \cite{ZakaDilauto}, \cite{ZakaPetersIso}, and from the properties of the construction of the reverse and inverse point, and sure from the Desargues condition, we have that,
\[
\begin{aligned}
-C'&=-\varphi(C)=\varphi(-C),\\
(A-C)'&=\varphi(A-C)=\varphi(A)-\varphi(C)=A'-C',\\
(B-C)'&=\varphi(B-C)=\varphi(B)-\varphi(C)=B'-C',\\
[(B-C)^{-1}]'&=\varphi[(B-C)^{-1}]=[\varphi(B-C)]^{-1}\\
&=[\varphi(B)-\varphi(C)]^{-1}. 
\end{aligned}
\]
If we consider, two 'three-vertex' $OP_2R$ and $O'P_2'\varphi(R)$, in Fig.\ref{Ratio3points.tr}, we see that they are \emph{Desarguesian}, for this reason, and for the results for construction of points $R'=([B-C]^{-1})'\cdot A'$ we have that,
\[
R'=\varphi(R)
\]
because, $P_2R||P_1(A-C)||P_1'(A-C)'=P_1'(A'-C') || P_2'\varphi(R)$, and parallelism is equivalence relation. So the points, $R'$ and $\varphi(R)$ are the cutting points of lines $\ell^{O'I'}$ and $\ell^{P_2'}_{P_1'(A-C)'}$, but this is a single point, for this, we have that $R'=\varphi(R).$ 

The case when, $\varphi: \ell^{OI} \to \ell^{OI}$ (there are traces $\ell^{OI}$), $\varphi$ can be seen as a composition of two translations $\varphi_1,\varphi_2$, with different trace of the line $\ell^{OI}$.
\qed

\begin{theorem}
The parallel projection between the two lines $\ell_1$ and $\ell_2$ in Desargues affine plane, \textbf{preserving the ratio} of $3-$points, 
\[
P_P(r(A,B;C))=r(P_P(A), P_P(B);P_P(C))
\]
\end{theorem}
\proof
If $\ell_1 || \ell_2$, we have that the parallel projection is a translation, and have true this theorem.\\
If lines $\ell_1$ and $\ell_2$  they are not parallel (so, they are cut at a single point), we have $A,B,C \in \ell_1$ and $P_P(A), P_P(B), P_P(C) \in \ell_2$. 

\begin{figure}[htbp]
	\centering
		\includegraphics[width=0.95\textwidth]{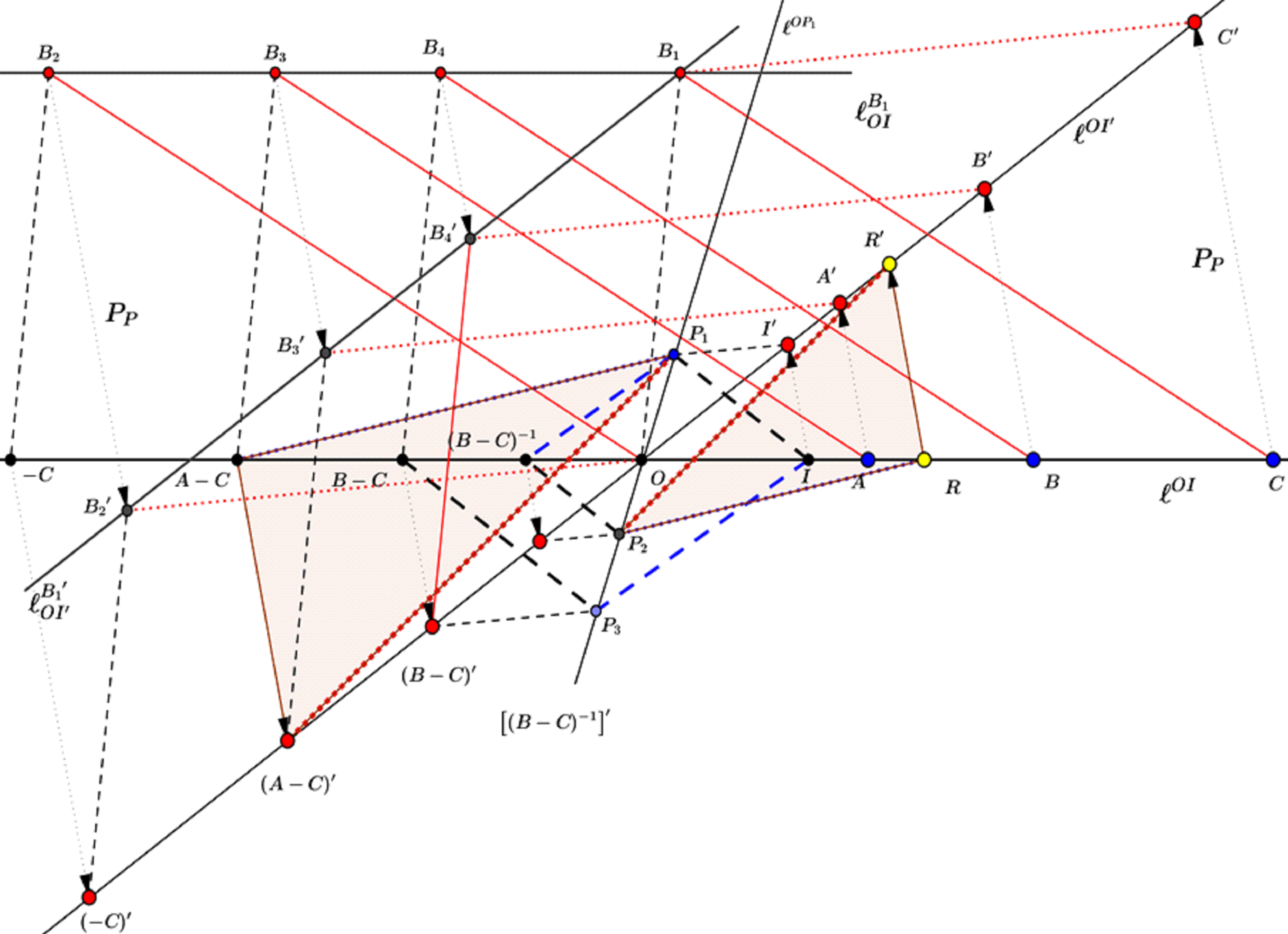}
	\caption{Ratio of 3-Points into parallel projection $P_P$ in Desargues affine plane.}
 \label{Ratio3points.pp}
\end{figure}

Lets have a parallel projection from the $\ell^{OI}-$line, to $\ell^{OI'}-$line (we are assuming that the lines intersect at the point $O$, if, the point of intersection would be another point, we can use a translation, to bring the point $O$ to the point of intersection) Fig.\ref{Ratio3points.pp}. 

So, have
\[ P_P: \ell^{OI} \to \ell^{OI'}.\] 
and $\forall A,B \in \ell^{OI}$, we have that $AP_P(A)\parallel BP_P(B)$, and $O'=P_P(O)=O$.

We mark with: $ I'=P_P(I)$, $A'=P_P(A)$, $B'=P_P(B)$, $C'=P_P(C)$, and  $(-C)'=P_P(-C)$,  $(A-C)'=P_P(A-C)$,  $(B-C)'=P_P(B-C)$, $\left[(B-C)^{-1}\right]'=P_P\left[(B-C)^{-1}\right]$.

Since addition of points do not depend on the position of the auxiliary point (see \cite{ZakaFilipi2016}, \cite{ZakaThesisPhd} ), we are keeping the same auxiliary point $B_1$, and construct the line $\ell^{B_1}_{OI'}$, and use the parallel projection with which, we get the points:
$P_P(B_1)=B_1$, $P_P(B_2)=B_2'$, $P_P(B_3)=B_3'$, $P_P(B_4)=B_4'$. Also, we know that the multiplication of points in a line of Desargues affine plane does not depend on the position of the auxiliary point (see \cite{ZakaThesisPhd}, \cite{FilipiZakaJusufi}), therefore for the multiplication of points in the line $\ell^{OI'}$, we we keep the same auxiliary points as in the case of multiplication of points, in the line $\ell^{OI}$, so the points $P_1, P_2, P_3$.

For parallel projection properties, we have the following parallelisms,
\[
\ell^{(-C)(-C)'} \parallel \ell^{(A-C)(A-C)'} \parallel \ell^{(B-C)(B-C)'} \parallel
\ell^{[(B-C)^{-1}][(B-C)^{-1}]'}  \parallel \ell^{II'} \]
continue
\[
 \ell^{II'} \parallel \ell^{AA'}  \parallel \ell^{RP_P(R)}  \parallel \ell^{BB'}  \parallel \ell^{CC'} \parallel \ell^{B_2B_2'} \parallel \ell^{B_3B_3'} \parallel \ell^{B_4B_4'}
\]
Since parallelism is an equivalence relation, we, from the above parallelisms, can distinguish
\[\ell^{(A-C)(A-C)'} \parallel \ell^{RP_P(R)}.
\]
For construction of ratio-point $R$, we have that, the point 
\[R=\ell^{OI}\cap \ell^{P_2}_{P_1(A-C)} \]
and we have the follow parallelism
\[ \ell^{(A-C)P_1} \parallel \ell^{P_2R}.\] 

Also, we have triads of collinear points $(A-C),O,R$ and $(A-C)',O,P_P(R)$, and $P_1,O,P_2$. Hence we have that the two three-vertexes 
\[
(A-C)P_1(A-C)'\quad \text{and}\quad RP_P(R)P_2\rightarrow\text{are Desarguesian.}
\]
For this we have that 
\[ \ell^{(A-C)'P_1} \parallel \ell^{P_2P_P(R)}. \]
It is easy to prove that, 
\[
\begin{aligned}
(-C)'&=P_P(-C)=-P_P(C), \\
(A-C)'&=P_P(A-C)=P_P(A)-P_P(C)=A'-C',\\
(B-C)'&=P_P(B-C)=P_P(B)-P_P(C)=B'-C',\\
\left[(B-C)^{-1}\right]'&=P_P\left[(B-C)^{-1}\right] \\
&=\left[P_P(B-C)\right]^{-1} \\
&=\left[P_P(B)-P_P(C)\right]^{-1} \\
&=\left[B'-C'\right]^{-1}.
\end{aligned}
\]

For construction of ratio point $R'$ in $\ell^{OI'}-$line, we have following parallelisms:
\[
\ell^{P_1[A'-C']}=\ell^{P_1(A-C)'}\parallel \ell^{P_2R'}
\]
and, from construction of point $R$, we had the parallelisms
\[
\ell^{P_1(A-C)}\parallel \ell^{P_2R}
\]
Also, we have triads of collinear points $(A-C),O,R$ and $(A-C)',O,R'$, and $P_1,O,P_2$. So, we have that the two three-vertexes,
\[
(A-C)P_1(A-C)'\quad \text{and}\quad RR'P_2\rightarrow\text{are Desarguesian.}
\]
For this we have that 
\[  \ell^{(A-C)'P_1} \parallel \ell^{P_2R'}.  \]
Also we have that the points $R'$ and $P_P(R)$ are in $\ell^{OI'}-$line.  So, we have that,
\[
\ell^{P_2P_P(R)} \parallel \ell^{P_2R}\quad\text{and points}\quad R',P_P(R)\in\ell^{OI'}
\]
Hence, have
\[R'=P_P(R).
\]
If, $P_P: \ell^{OI} \to \ell^{OI}$ we have the \emph{identical-parallel projection}. The case where $P_P: \ell^{OI} \to \ell^{O'I'}$ , and $\ell^{OI} \parallel \ell^{O'I'}$, the parallel projection is a translation.
\qed
\begin{remark}
If lines $\ell_1 = \ell_2$, then the parallel projection is a translation, or identical transform. 
\end{remark}

\begin{theorem}
dilatation $\delta$ preserve the ratio of this points, so
\[
\delta(r(A,B;C))=r(\delta(A), \delta(B);\delta(C))
\]
\end{theorem}
\proof
Lets have firstly, a dilatation with an fixed point $V \notin \ell^{OI}$ in Desargues affine plane ($\delta(V)=V$), which $\delta: \ell^{OI} \to \ell^{O'I'}$ we know that dilatation preserves parallelism, therefore $\ell^{OI} \parallel \ell^{O'I'}$.

We mark with: $O'=\delta(O)$, $I'=\delta(I)$, $A'=\delta(A)$, $B'=\delta(B)$, $C'=\delta(C)$, $\delta(R)$. Also, dilatation $\delta$, we apply it to auxiliary points, 
$B_1'=\delta(B_1)$, $B_2'=\delta(B_2)$, $B_3'=\delta(B_3)$, $B_4'=\delta(B_4)$ and mark $P_1'=\delta(P_1)$, $P_2'=\delta(P_2)$, $P_3'=\delta(P_3)$. With results in \cite{ZakaDilauto}, \cite{ZakaPetersIso}, and  from the properties of the construction of the inverse point, we have that,
\[ \begin{aligned}
(-C)'&=\delta(-C)=-\delta(C)=-C'  \\
(A-C)'&=\delta(A-C)=\delta(A)-\delta(C)=A'-C'\\
(B-C)'&=\delta(B-C)=\delta(B)-\delta(C)=B'-C'\\
[(B-C)^{-1}]'&=\delta[(B-C)^{-1}]=[\delta(B-C)]^{-1}
=[\delta(B)-\delta(C)]^{-1}=[B'-C']^{-1}
\end{aligned}
\]
Now construct the ratio point $R'=r(A',B';C')$ in $\delta(\ell^{OI})-$line. From definition of ratio of three points and the above notes, we have that,
\[R'=(B'-C')^{-1}(A'-C')
\]
is ratio point in $\delta(\ell^{OI})-$line.

\begin{figure}[htbp]
	\centering
		\includegraphics[width=1 \textwidth]{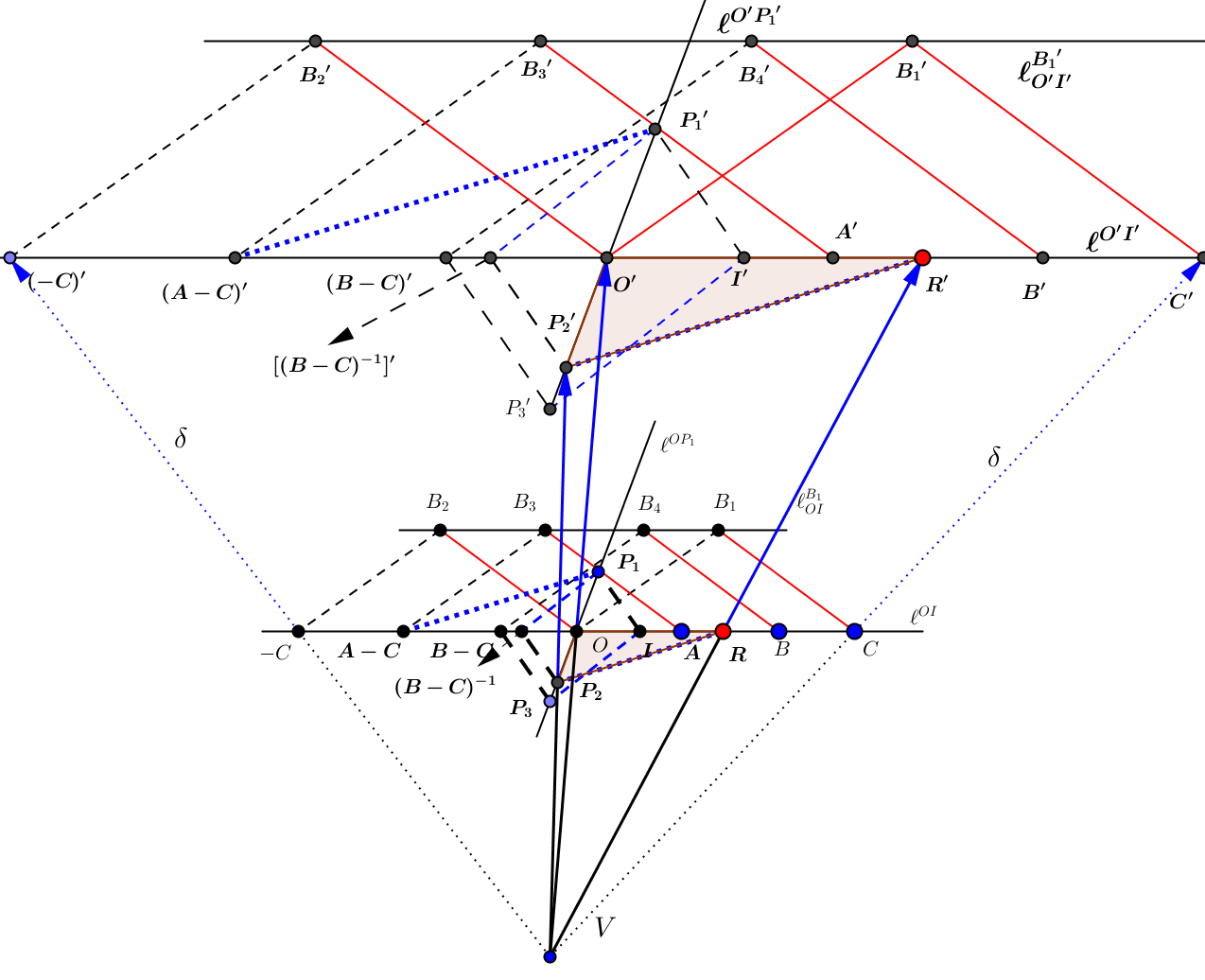}
	\caption{Ratio of 3-Points into a dilatation $\delta$ with fixed point $V\notin \ell^{OI}-$line, in Desargues affine plane.}
		\label{Ratio3points.dil}
\end{figure}

If we consider, two 'three-vertex' $OP_2R$ and $O'P_2'\delta(R)$, in Fig.\ref{Ratio3points.dil}, we see that they are Desarguesian, for this reason, and for the results for construction of points $R'=(B'-C')^{-1}(A'-C')$ we have that,
\[
R'=\delta(R)
\]
because, $P_2R||P_1(A-C)||P_1'(A-C)'=P_1'(A'-C') || P_2'R'$,
also from dilatation properties we have that 
\[P_2R||\delta(P_2)\delta(R)=P_2'\delta(R),\]
and parallelism is equivalence relation, hence we have 
\[P_2'R' \parallel P_2'\delta(R),\]
and the points $R', \delta(R) \in \delta(\ell^{OI})=\ell^{O'I'}$.

So the points, $R'$ and $\delta(R)$ are the cutting points of lines $\ell^{O'I'}$ and $\ell^{P_2'}_{P_1'(A-C)'}$ (the line which passes from point $P_2'$ and is parallel with line $\ell^{P_1'(A-C)'}$), so

\[R'=\ell^{O'I'} \cap \ell^{P_2'}_{P_1'(A-C)'} \quad 
\text{and} \quad \delta(R)=\ell^{O'I'} \cap \ell^{P_2'}_{P_1'(A-C)'} \]

thus, this is a single point, for this, we have that $R'=\delta(R)$. Hence
\[
\delta[r(A,B;C)]=r[\delta(A), \delta(B);\delta(C)].
\]

\begin{figure}[htbp]
	\centering
		\includegraphics[width=1 \textwidth]{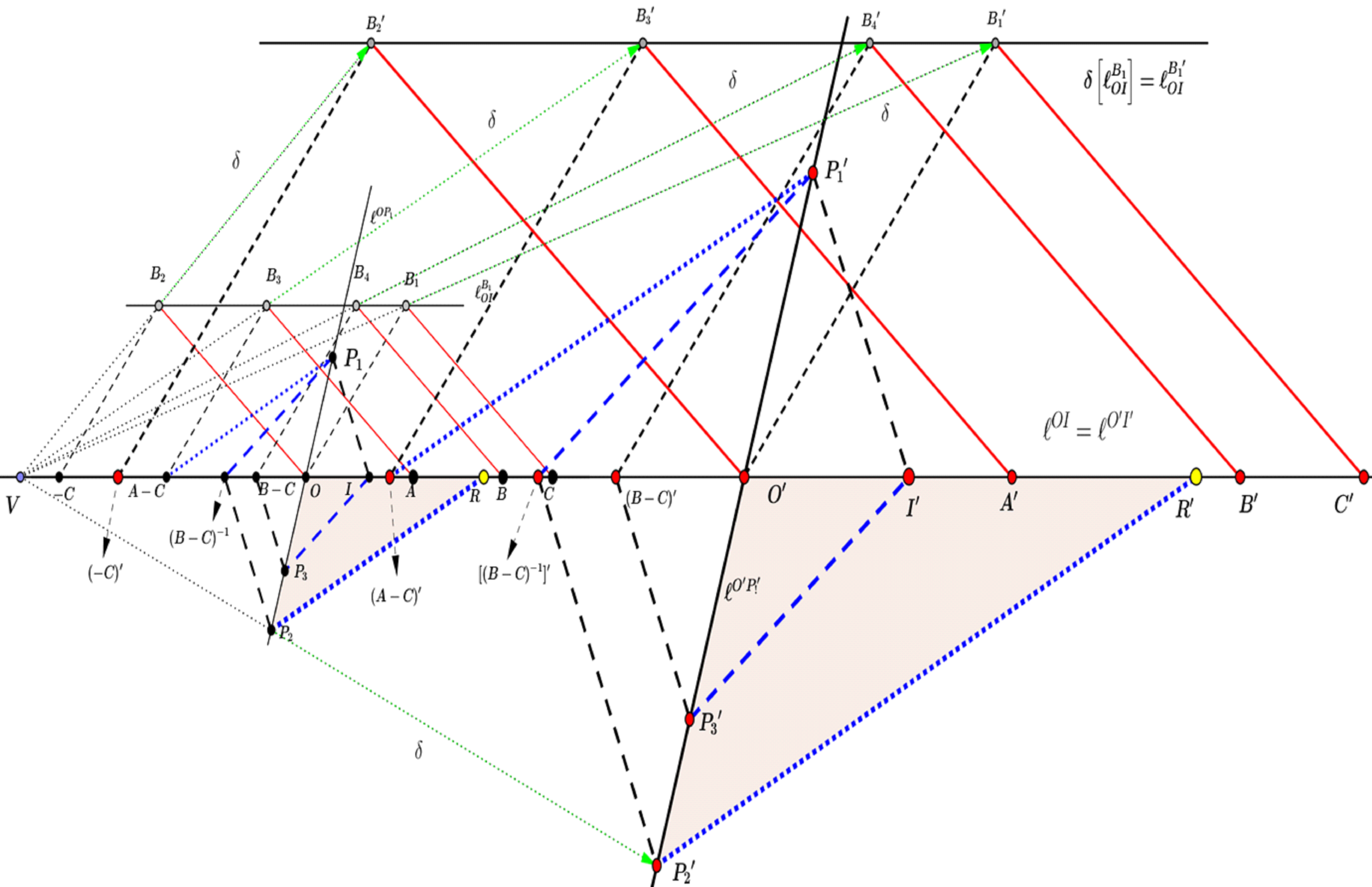}
	\caption{Ratio of 3-Points into a dilatation $\delta$ with fixed point $V\in \ell^{OI}-$line, in Desargues affine plane.}
	\label{ratio3points.dil2}
\end{figure}

Let's now consider the case, when $\delta: \ell^{OI} \to \ell^{OI}$ (is case where the fixed point $V \in \ell^{OI}$), Fig.\ref{ratio3points.dil2}. 

We mark (same as in the first case) with: $O'=\delta(O)$, $I'=\delta(I)$, $A'=\delta(A)$, $B'=\delta(B)$, $C'=\delta(C)$, $\delta(R)$, all this points are in $\ell^{OI}-$line. Also, dilatation $\delta$, we apply it to auxiliary points, 
$B_1'=\delta(B_1)$, $B_2'=\delta(B_2)$, $B_3'=\delta(B_3)$, $B_4'=\delta(B_4)$ (from dilatation properties have that, $\ell^{B_1}_{OI}\parallel \delta \left[  \ell^{B_1}_{OI} \right]=\ell^{B_1'}_{OI}$),  also marked $P_1'=\delta(P_1)$, $P_2'=\delta(P_2)$, $P_3'=\delta(P_3)$. With results in \cite{ZakaDilauto}, \cite{ZakaPetersIso}, and  from the properties of the construction of the inverse point, we have that,
\[ \begin{aligned}
(-C)'&=\delta(-C)=-\delta(C)=-C'  \\
(A-C)'&=\delta(A-C)=\delta(A)-\delta(C)=A'-C'\\
(B-C)'&=\delta(B-C)=\delta(B)-\delta(C)=B'-C'\\
[(B-C)^{-1}]'&=\delta[(B-C)^{-1}]=[\delta(B-C)]^{-1}
=[\delta(B)-\delta(C)]^{-1}=[B'-C']^{-1}
\end{aligned}
\]
Now construct the ratio point $R'=r(A',B';C')$ in $\delta(\ell^{OI})\equiv \ell^{OI}-$line. From definition of ratio of three points and the above notes, we have that,
\[R'=(B'-C')^{-1}(A'-C')
\]
is ratio point in $\ell^{OI}=\delta(\ell^{OI})-$line.

If we consider, two 'three-vertex' $OP_2R$ and $O'P_2'\delta(R)$, in Fig.\ref{Ratio3points.dil}, we see that they are Desarguesian, for this reason, and for the results for construction of points $R'=(B'-C')^{-1}(A'-C')$ we have that,
\[
R'=\delta(R) 
\]
because, $P_2R||P_1(A-C)||\delta(P_1)\delta(A-C)=P_1'(A-C)'=P_1'(A'-C') || P_2'R'$, 
also from dilatation properties we have that 
\[ P_2R||\delta(P_2)\delta(R)=P_2'\delta(R), \]
and parallelism is equivalence relation, hence we have 
\[ P_2'R' \parallel P_2'\delta(R),\] 
and we have that the points $R', \delta(R) \in \delta(\ell^{OI})=\ell^{OI}$.

So the points, $R'$ and $\delta(R)$ are the cutting points of lines $\ell^{OI}$ and $\ell^{P_2'}_{P_1'(A-C)'}$ (the line which passes from point $P_2'$ and is parallel with line $\ell^{P_1'(A-C)'}$), so

\[R'=\ell^{O'I'} \cap \ell^{P_2'}_{P_1'(A-C)'} \quad 
\text{and} \quad \delta(R)=\ell^{O'I'} \cap \ell^{P_2'}_{P_1'(A-C)'} \]

thus, this is a single point, for this, we have that $R'=\delta(R)$. Hence
\[
\delta[r(A,B;C)]=r[\delta(A), \delta(B);\delta(C)].
\]
\qed

\bibliographystyle{amsplain}
\bibliography{RCRrefs}

\end{document}